\documentclass[a4paper,11pt]{article}
\usepackage{amsmath, amsfonts, amscd, amssymb, amsthm, enumerate, xypic}

\def\Diag{\text{Diag}}
\def\Mat{\text{M}}
\def\NT{\text{NT}}
\def\GL{\text{GL}}

\def\id{\text{id}}

\newcommand{\Ker}{\operatorname{Ker}}
\newcommand{\Vect}{\operatorname{Span}}
\newcommand{\im}{\operatorname{Im}}

\newcommand{\Sp}{\operatorname{Sp}}
\newcommand{\codim}{\operatorname{codim}}
\renewcommand{\setminus}{\smallsetminus}

\def\K{\mathbb{K}}
\def\R{\mathbb{R}}
\def\C{\mathbb{C}}

\def\Z{\mathbb{Z}}

\def\DG{\text{DG}}

\def\calS{\mathcal{S}}

\def\lcro{\mathopen{[\![}}
\def\rcro{\mathclose{]\!]}}

\theoremstyle{definition}
\newtheorem{Def}{Definition}

\theoremstyle{plain}
\newtheorem{theo}{Theorem}
\newtheorem{prop}[theo]{Proposition}

\newtheorem{lemme}[theo]{Lemma}
\newtheorem{claim}{Claim}

\theoremstyle{plain}

\theoremstyle{remark}
\newtheorem{Rems}{Remarks}
\newtheorem{Rem}[Rems]{Remark}

\title{The linear preservers of real diagonalizable matrices}
\author{Bernard Rand\'e\footnote{Professor of Mathematics at Lyc\'ee Louis le Grand, 123, rue Saint-Jacques, 75005 Paris, FRANCE. e-mail: bernard.rande@free.fr} \hskip 2mm and \hskip 2mm
Cl\'ement de Seguins Pazzis\footnote{Professor of Mathematics at Lyc\'ee Priv\'e Sainte-Genevi\`eve, 2, rue
de l'\'Ecole des Postes, 78029 Versailles Cedex, FRANCE. e-mail: dsp.prof@gmail.com}
}

\begin{document}

\thispagestyle{plain}
\maketitle

\begin{abstract}
Using a recent result of Bogdanov and Guterman on the linear preservers of pairs of simultaneously diagonalizable matrices, we
determine all the automorphisms of the vector space $\Mat_n(\R)$ which stabilize the set of diagonalizable matrices.
To do so, we investigate the structure of linear subspaces of diagonalizable matrices of $\Mat_n(\R)$
with maximal dimension.
\end{abstract}



\section{Introduction}

In this paper, $\K$ will denote an arbitrary field and $\K^*:=\K \setminus \{0\}$.
Given positive integers $n$ and $p$, we let $\Mat_{n,p}(\K)$ denote the set of matrices with $n$ rows, $p$ columns and entries in $\K$.
Given a positive integer $n$, we let $\Mat_n(\K)$ denote the vector space of
square matrices of order $n$ with entries in $\K$; we let $D_n(\K)$ denote its linear subspace of \emph{diagonal} matrices,
$\DG_n(\K)$ the subset of \emph{diagonalizable} matrices, $\NT_n(\K)$ the linear subspace of strictly upper triangular matrices,
$\calS_n(\K)$ the linear subspace of symmetric matrices, and $O_n(\K)$ the subgroup of matrices $A$ for which $A^TA=I_n$.
We will systematically use the basic fact that $O\calS_n(\K) O^{-1}=\calS_n(\K)$ for every $O \in O_n(\K)$.
For $(i,j)\in \lcro 1,n\rcro^2$, we let $E_{i,j}$ denote the elementary matrix of $\Mat_n(\K)$
with entry zero everywhere except in position $(i,j)$ where the entry is $1$. Given a square matrix
$A \in \Mat_n(\K)$, we let $\Sp(A)$ denote its \textbf{spectrum}, i.e.\ the set of eigenvalues of $A$ in $\K$,
and for $\lambda \in \Sp(A)$, we let $E_\lambda(A):=\Ker(A-\lambda\cdot I_n)$ denote the eigenspace associated to $\lambda$ and $A$.

Linear preserver problems have been a very active field of research in the recent decades.
Although the first results date back to Frobenius in the 19th century
(see \cite{Frobenius} on the linear preservers of the determinant) and Dieudonn\'e around 1950
(see \cite{Dieudonne} on the linear preservers of non-singularity), most of the known linear preserver theorems have
been established in the last thirty years.
Amongst the presently unresolved issues is the determination of the linear preservers of diagonalizability over an arbitrary field:
is there a neat explicit description of the automorphisms $f$ of the vector space $\Mat_n(\K)$
which stabilize $\DG_n(\K)$ i.e.\ such that $f(M)$ is diagonalizable for every diagonalizable $M \in \Mat_n(\K)$?
Given a non-singular matrix $P \in \GL_n(\K)$, a linear form $\lambda$ on $\Mat_n(\K)$ and a scalar $\mu \in \K$,
the maps
$$\varphi_{\lambda,P,\mu} : M \mapsto \lambda(M)\cdot I_n+\mu\,PMP^{-1} \quad \text{and} \quad
\psi_{\lambda,P,\mu} : M \mapsto \lambda(M)\cdot I_n+\mu\,PM^TP^{-1}.$$
are both endomorphisms of the vector space $\Mat_n(\K)$ which preserve diagonalizability,
and they are automorphisms if and only if $\mu \neq 0$ and $\lambda(I_n) \neq -\mu$; the obvious conjecture is that
those are the only automorphisms of the vector space $\Mat_n(\K)$
which preserve diagonalizability. So far, this has only been established for algebraically closed fields
of characteristic $0$ (\cite{Omladic}) using the Motzkin-Taussky theorem (\cite{MoTau1} and \cite{MoTau2}).

Contrary to classical linear preserver theorems such as the ones mentioned earlier or the Botta-Pierce-Watkins theorem on the linear preservers of nilpotency, it is reasonable to think that the ground field plays an important role in the structure of the preservers of diagonalizability
since it already has a major impact on the geometry of the set of diagonalizable matrices: for example
the subspace $\calS_n(\R)$ of symmetric matrices of $\Mat_n(\R)$ only consists of diagonalizable matrices
and has dimension $\dbinom{n+1}{2}$ whilst $\Mat_n(\C)$ does not contain any subspace having this property if $n \geq 2$.
On the contrary, the Motzkin-Taussky theorem states that any linear subspace of diagonalizable matrices of $\Mat_n(\C)$
has simultaneously diagonalizable elements hence is conjugate to a linear subspace of $D_n(\C)$.

\subsection{Main results}

Here is our main theorem:

\begin{theo}\label{linprestheo}
Let $n \geq 2$ be an integer and assume that every symmetric matrix of $\Mat_n(\K)$ is diagonalizable.
Let $f : \Mat_n(\K) \rightarrow \Mat_n(\K)$ be a vector space automorphism which stabilizes $\DG_n(\K)$.
Then there exists a non-singular matrix $P \in \GL_n(\K)$, a linear form $\lambda : \Mat_n(\K) \rightarrow \K$
and a scalar $\mu \in \K^*$
such that $\lambda(I_n) \neq -\mu$ and $f=\varphi_{\lambda,P,\mu}$ or $f=\psi_{\lambda,P,\mu}$.
\end{theo}

Notice in particular that the above theorem holds for $\K=\R$ and more generally for every
intersection of real closed fields (it is known that those fields are precisely the ones for which
every symmetric matrix of any order is diagonalizable \cite{Waterhouse}).

In order to prove Theorem \ref{linprestheo}, we will study large linear subspaces of $\Mat_n(\K)$ consisting only of
diagonalizable matrices.

\begin{Def}
A linear subspace $V$ of $\Mat_n(\K)$ will be called \textbf{diagonalizable} if all its elements are diagonalizable.
\end{Def}

Beware that this \emph{does not mean} that the matrices of $V$ are \textbf{simultaneously diagonalizable},
i.e.\ that $PVP^{-1} \subset D_n(\K)$ for some $P \in \GL_n(\K)$.
Given a diagonalizable subspace $V$ of $\Mat_n(\K)$, we have $V \cap \NT_n(\K)=\{0\}$
since only $0$ is both diagonalizable and nilpotent. It follows that $\dim V \leq \codim_{\Mat_n(\K)} \NT_n(\K)$.
This yields:

\begin{prop}
Let $V$ be a diagonalizable subspace of $\Mat_n(\K)$. Then \\
$$\dim V \leq \dbinom{n+1}{2}.$$
\end{prop}

Notice that this upper bound is tight for $\K=\R$ since $S_n(\R)$ is a diagonalizable subspace of $\Mat_n(\R)$
with dimension $\dbinom{n+1}{2}$. However, for other fields, the previous upper bound may not be reached.
For example, if $\K$ is algebraically closed of characteristic $0$,
an easy consequence of the Motzkin-Taussky theorem (\cite{MoTau1,MoTau2})
is that the largest dimension of a diagonalizable subspace of $\Mat_n(\K)$ is $n$.

\begin{Def}
A diagonalizable subspace of $\Mat_n(\K)$ will be called \textbf{maximal} if its dimension is $\dbinom{n+1}{2}$.
\end{Def}

Note that this definition of maximality is quite different from the notion of maximality with respect to the inclusion of diagonalizable
subspace. We will discuss this in Section \ref{last}.
Recall that the field $\K$ is called:
\begin{itemize}
\item \textbf{formally real} if $-1$ is not a sum of squares in $\K$;
\item \textbf{Pythagorean} if any sum of two squares is a square in $\K$.
\end{itemize}
With that in mind, here is our major result on maximal diagonalizable subspaces.

\begin{theo}\label{structmaxsub}
Let $n \geq 2$. Assume there exists a maximal diagonalizable subspace of $\Mat_n(\K)$.
Then:
\begin{enumerate}[(a)]
\item the field $\K$ is formally real and Pythagorean;
\item every maximal diagonalizable subspace $V$ of $\Mat_n(\K)$ is \textbf{conjugate} to $\calS_n(\K)$,
i.e.\ there exists a $P \in \GL_n(\K)$ such that
$V=P\,\calS_n(\K)\,P^{-1}$;
\item every symmetric matrix of $\Mat_n(\K)$ is diagonalizable.
\end{enumerate}
\end{theo}

In particular, if $\Mat_n(\K)$ contains a maximal diagonalizable subspace, then $\calS_n(\K)$ is automatically such a subspace
and those subspaces make a single orbit under conjugation.

\subsection{Structure of the paper}

In Section \ref{linsub}, Theorem \ref{structmaxsub} will be proved by induction on $n$.
We will then derive Theorem \ref{linprestheo} (Section \ref{linpres}) by appealing to a recent result of Bogdanov and Guterman,
and will also investigate there the strong linear preservers of diagonalizability.

\section{Linear subspaces of diagonalizable matrices}\label{linsub}

In order to prove Theorem \ref{structmaxsub}, we will proceed by induction: the case
$n=2$ will be dealt with in Section \ref{n=2section} and the rest of the induction will be carried out in Section \ref{ngeq3section}.

\subsection{The case $n=2$}\label{n=2section}

Let $V$ be a diagonalizable subspace of $\Mat_2(\K)$ with dimension $3$.
We first wish to prove that $V$ is conjugate to $\calS_2(\K)$. \\
Notice first that $I_2 \in V$: indeed, if $I_2 \not\in V$, then
$\Mat_2(\K)=\Vect(I_2)\oplus V$ hence every matrix of $\Mat_2(\K)$ would be diagonalizable, which of course is not the case for
$\begin{bmatrix}
0 & 1 \\
0 & 0
\end{bmatrix}$. Choose then some $A \in V \setminus \Vect(I_2)$. Since $A$ is diagonalizable, we lose
no generality assuming that $A$ is diagonal (we simply conjugate $V$ with an appropriate matrix of $\GL_2(\K)$).
In this case, $(I_2,A)$ is a basis of $D_2(\K)$ hence $D_2(\K) \subset V$.
Choose some $B \in V \setminus D_2(\K)$ and write it
$B=\begin{bmatrix}
a & b \\
c & d
\end{bmatrix}$. Replacing $B$ with $B-\begin{bmatrix}
a & 0 \\
0 & d
\end{bmatrix}$, we lose no generality assuming $B=\begin{bmatrix}
0 & b \\
c & 0
\end{bmatrix}$. Then $B$ is diagonalizable and non-zero hence $b \neq 0$ and $c \neq 0$.
Multiplying $B$ with an appropriate scalar, we may then assume $c=1$.
However, $B$ is diagonalizable and has trace $0$ hence its eigenvalues are $\lambda$ and $-\lambda$ for some $\lambda\in \K^*$:
we deduce that $b=-\det B=\lambda^2$. Conjugating $V$ with the diagonal matrix $\begin{bmatrix}
1 & 0 \\
0 & \lambda
\end{bmatrix}$, we are reduced to the case $V$ contains the subspace $D_2(\K)$ and the matrix
$\begin{bmatrix}
0 & \lambda \\
\lambda & 0
\end{bmatrix}$: in this case, we deduce that $V$ contains $\calS_2(\K)$ and equality of dimensions shows that $V=\calS_2(\K)$, QED.

\vskip 2mm
\noindent We will conclude the case $n=2$ by reproving the following classical result:

\begin{prop}
The subspace $\calS_2(\K)$ of $\Mat_2(\K)$ is diagonalizable if and only if $\K$ is formally real and Pythagorean.
\end{prop}

\begin{proof}
Assume $\calS_2(\K)$ is diagonalizable. Let $(a,b)\in \K^2$. The matrix
$\begin{bmatrix}
a & b \\
b & -a
\end{bmatrix}$ is then diagonalizable with trace $0$: its eigenvalues are then $\lambda$ and $-\lambda$ for some $\lambda \in \K$.
Computing the determinant yields $a^2+b^2=\lambda^2$. Moreover, if $a^2+b^2=0$, then $\lambda=0$
hence $\begin{bmatrix}
a & b \\
b & -a
\end{bmatrix}=0$ which shows that $a=b=0$. This proves that $\K$ is formally real and Pythagorean. \\
Conversely, assume $\K$ is formally real and Pythagorean. Then $\K$ has characteristic $0$ hence we may split
any matrix of $\calS_2(\K)$ as $c\cdot I_2+\begin{bmatrix}
a & b \\
b & -a
\end{bmatrix}$ for some $(a,b,c)\in \K^3$. Fixing an arbitrary $(a,b)\in \K^2 \setminus \{(0,0)\}$,
it will suffice to prove that $A:=\begin{bmatrix}
a & b \\
b & -a
\end{bmatrix}$ is diagonalizable. However, its characteristic polynomial is $\chi_A=X^2-(a^2+b^2)$
and $a^2+b^2=d^2$ for some $d \in \K^*$ since $\K$ is formally real and Pythagorean.
Hence $\chi_A=(X-d)(X+d)$ with $d \neq -d$, which proves that $A$ is diagonalizable.
\end{proof}

\subsection{The case $n \geq 3$}\label{ngeq3section}

Our proof will feature an induction on $n$ and a reduction to special cases.
First, we will introduce a few notations. Given a maximal diagonalizable subspace $V$ of $\Mat_n(\K)$,
we write every matrix of it as
\begin{multline*}
M=\begin{bmatrix}
K(M) & C(M) \\
L(M) & \alpha(M)
\end{bmatrix} \quad \text{with $K(M) \in \Mat_{n-1}(\K)$, $L(M) \in \Mat_{1,n-1}(\K)$,} \\
\text{$C(M) \in \Mat_{n-1,1}(\K)$ and $\alpha(M) \in \K$.}
\end{multline*}

\subsubsection{The fundamental lemma}

\begin{lemme}\label{imagelemma}
Let $p \in \lcro 1,n-1\rcro$ and let $A \in \Mat_p(\K)$, $B \in \Mat_{n-p}(\K)$ and
$C \in \Mat_{p,n-p}(\K)$. Assume that
$M=\begin{bmatrix}
A & C \\
0 & B
\end{bmatrix}$ is diagonalizable. Then, for every eigenvalue $\lambda$ of $B$ and every corresponding eigenvector $X$ of $B$,
one has $CX\in \im(A-\lambda\cdot I_p)=\underset{\mu \in \Sp(A)\setminus\{\lambda\}}{\bigoplus}\Ker(A-\mu\cdot I_p)$.
\end{lemme}

\begin{proof}
We lose no generality considering only the case $\lambda=0$.
Set then $F:=\K^p \times \{0\}$ seen as a linear subspace of $\K^n$
and let $u$ denote the endomorphism of $\K^n$ canonically associated to $M$.
Let $X \in \Ker B$. Then $x:=\begin{bmatrix}
0 \\
X
\end{bmatrix}$ satisfies $u(x)=\begin{bmatrix}
CX \\
0
\end{bmatrix}$ and we wish to prove that $u(x) \in \im u_{|F}$.
It thus suffices to prove that $F \cap \im u \subset \im u_{|F}$.
However, $u_{|F}$ is diagonalizable, so we may choose a basis $(e_1,\dots,e_p)$ of $F$
consisting of eigenvectors of $u$, and extend it to a basis $(e_1,\dots,e_n)$ of $\K^n$
consisting of eigenvectors of $u$. The equality $\im u_{|F}=F \cap \im u$
follows easily by writing $\im u=\Vect\{e_i \mid 1 \leq i \leq n \; \text{such that} \; u(e_i) \neq 0\}$
and $\im u_{|F}=\Vect\{e_i \mid 1 \leq i \leq p \; \text{such that} \; u(e_i) \neq 0\}$.
\end{proof}

\begin{Rem}
The previous lemma may also be seen as an easy consequence of Roth's theorem (see \cite{Roth}).
\end{Rem}

\subsubsection{A special case}

\begin{prop}\label{finalisation}
Assume $\K$ is formally real and Pythagorean.
Let $V$ be a maximal diagonalizable subspace of $\Mat_n(\K)$. Assume furthermore that:
\begin{enumerate}[(i)]
\item every matrix $M \in V$ with zero as last row
has the form $\begin{bmatrix}
S & ? \\
0 & 0
\end{bmatrix}$ for some symmetric matrix $S \in \calS_{n-1}(\K)$;
\item for every symmetric matrix $S \in \calS_{n-1}(\K)$,
the subspace $V$ contains a unique matrix of the form
$\begin{bmatrix}
S & ? \\
0 & 0
\end{bmatrix}$;
\item the subspace $V$ contains $E_{n,n}$.
\end{enumerate}
Then $V$ is conjugate to $\calS_n(\K)$.
\end{prop}

We start with a first claim:

\begin{claim}
The subspace $V$ contains $\begin{bmatrix}
S & 0 \\
0 & 0
\end{bmatrix}$ for any $S \in \calS_{n-1}(\K)$.
\end{claim}

Let indeed $S \in \calS_{n-1}(\K)$. We already know that there is a column $C_1 \in \Mat_{n-1,1}(\K)$ such that
$V$ contains $\begin{bmatrix}
S & C_1 \\
0 & 0
\end{bmatrix}$. It follows that $S$ is diagonalizable. Since $E_{n,n} \in V$, we know that $\begin{bmatrix}
S & C_1 \\
0 & \lambda
\end{bmatrix}$ is in $V$ and is thus diagonalizable for every $\lambda \in \Sp(S)$.
Lemma \ref{imagelemma} then shows that $C_1 \in \im(S-\lambda\cdot I_n)$ for every $\lambda \in \Sp(S)$.
Since $S$ is diagonalizable, we have $\underset{\lambda \in \Sp(S)}{\bigcap}\,\im(S-\lambda\cdot I_n)=\{0\}$
hence $C_1=0$. This proves our first claim.

\vskip 4mm
Moreover, assumptions (i) and (ii) show that
the kernel of $L : M \mapsto L(M)$ has dimension at most $1+\dim \calS_{n-1}(\K)=\binom{n+1}{2}-(n-1)$,
and it follows from the rank theorem that $L(V)=\Mat_{1,n-1}(\K)$.
Notice also that $S_{n-1}(\K) \oplus \NT_{n-1}(\K)=\Mat_{n-1}(\K)$.
Using this, we find a linear endomorphism $u$ of $\Mat_{1,n-1}(\K)$
and a linear map $v : \Mat_{1,n-1}(\K) \rightarrow \NT_{n-1}(\K)$
such that
$V$ contains
$M_L=\begin{bmatrix}
v(L) & u(L)^T \\
L & 0
\end{bmatrix}$ for every $L \in \Mat_{1,n-1}(\K)$.
Our next claim follows:

\begin{claim}
There is a scalar $\lambda \in \K^*$ such that $u=\lambda^2\cdot \id$,
and $v=0$.
\end{claim}

We will show indeed that $u$ is diagonalizable with a single eigenvalue which is a square.
We start by considering the row matrix $L_1=\begin{bmatrix}
0 & \cdots & 0 & 1
\end{bmatrix} \in \Mat_{1,n-1}(\K)$. We may then write
$$M_L=\begin{bmatrix}
U_1 & C_1 & C_2 \\
0 & 0 & a \\
0 & 1 & 0
\end{bmatrix} \quad \text{for some $U_1 \in \NT_{n-2}(\K)$, some $a \in K$, and some $(C_1,C_2)\in \Mat_{n-2,1}(\K)^2$.}$$
The matrix $U_1$ is then both diagonalizable and nilpotent hence $U_1=0$.
The matrix $B:=\begin{bmatrix}
0 & a \\
1 & 0
\end{bmatrix}$ must also be diagonalizable hence $a=\lambda^2$ for some $\lambda \in \K^*$ (see the proof of Section \ref{n=2section}).
This shows that the eigenvalues of $B$ are $\lambda$ and $-\lambda$.
Consider the matrix $C=\begin{bmatrix}
C_1 & C_2 \\
\end{bmatrix}$. Choose an eigenvector $X$ of $B$ which corresponds to the eigenvalue $\lambda$. The matrix
$\begin{bmatrix}
\lambda\cdot I_{n-2} & C \\
0 & B
\end{bmatrix}=\begin{bmatrix}
\lambda\cdot I_{n-2} & 0 \\
0 & 0
\end{bmatrix}+M_L$ belongs to $V$ and is thus diagonalizable. Using Lemma \ref{imagelemma}, we find that
$CX=0$. Since the eigenvectors of $B$ span $\K^2$, we deduce that $C=0$.
This proves that $v(L)=0$ and $L_1$ is an eigenvector of $u$ for the eigenvalue $\alpha=\lambda^2$.

\vskip 2mm
We may now prove that any non-zero vector of $\Mat_{1,n-1}(\K)$ is an eigenvector of $u$. \\
We identify canonically $\Mat_{1,n-1}(\K)$ with $\K^{n-1}$ and equip it with the canonical regular quadratic form
$q : (x_1,\dots,x_{n-1}) \mapsto \sum_{k=1}^{n-1} x_k^2$.
The Witt theorem (see Theorem XV 10.2 in \cite{Lang}) then shows that $O_{n-1}(\K)$ acts transitively on the \emph{sphere} $q^{-1}\{1\}$.
Let $L \in \Mat_{1,n-1}(\K)$ be such that $q(L)=1$. Then there exists some $O \in O_{n-1}(\K)$ such that
$L=L_1O$. For the non-singular matrix
$P:=\begin{bmatrix}
O & 0 \\
0 & 1
\end{bmatrix}$, the subspace $PVP^{-1}$ satisfies all the assumptions from Proposition \ref{finalisation}
and it contains the matrix
$\begin{bmatrix}
Ov(L)O^T & O u(L)^T \\
L_1 & 0
\end{bmatrix}$. From the previous step, we deduce that $O v(L)O^T$ is symmetric, hence $v(L)$ also is, which proves that $v(L)=0$.
Moreover, we find that $O u(L)^T=\beta L_1^T$ for some $\beta \in \K$ hence
$u(L)=\beta L_1O=\beta L$. We deduce that every vector of the sphere $q^{-1}\{1\}$ is an eigenvector of
$u$. However, since $\K$ is Pythagorean and formally real, we find that for every non-zero vector $x \in \Mat_{1,n-1}(\K)$,
there is a scalar $\mu \neq 0$ such that $q(x)=\mu^2$ hence $\frac{1}{\mu}\,x$ is an eigenvector of $u$
in $\Ker v$ and so is $x$. This shows that $v=0$ and $u$ is a scalar multiple of the identity,
hence $u=\lambda^2\cdot \id$ by the above notations.

\vskip 3mm
\noindent We may now conclude. With the previous $\lambda$, set $P:=\begin{bmatrix}
I_{n-1} & 0 \\
0 & \lambda
\end{bmatrix}$ and notice that the subspace $V_1=PVP^{-1}$ satisfies all the conditions from Proposition \ref{finalisation}
with the additional one:
\begin{center}
for every $L \in \Mat_{1,n-1}(\K)$, the subspace $V_1$ contains $\begin{bmatrix}
0 & L^T \\
L & 0
\end{bmatrix}$.
\end{center}
It easily follows that $\calS_n(\K) \subset V_1$ hence $\calS_n(\K)=V_1$ since $\dim V_1=\dbinom{n+1}{2}=\dim \calS_n(\K)$.
This finishes the proof of Proposition \ref{finalisation}.
\hfill $\square$

\subsubsection{Proof of Theorem \ref{structmaxsub} by induction}

We now proceed by induction. Let $n \geq 3$ and assume that
if there is a maximal diagonalizable subspace of $\Mat_{n-1}(\K)$, then this subspace is conjugate to $\calS_{n-1}(\K)$
and $\K$ is formally real and Pythagorean (note that \emph{we do not assume that such a subspace exists at this point}).
Let $V$ be a maximal diagonalizable subspace of $\Mat_n(\K)$.

Consider the subspace $W=\Ker L$ consisting of the matrices $M$ of $V$ such that $L(M)=0$.
In $W$, consider the subspace $W'$ of matrices $M$ such that $K(M)=0$.
Using the rank theorem twice shows that:
$$\dim V=\dim L(V)+\dim K(W)+\dim W'.$$
\begin{itemize}
\item We readily have $\dim L(V) \leq n-1$.
\item Notice that $M \mapsto \alpha(M)$ is injective on $W'$: indeed, any matrix $M \in W'$ such that $\alpha(M)=0$
has the form $M=\begin{bmatrix}
0 & C(M) \\
0 & 0
\end{bmatrix}$, hence is nilpotent, but also diagonalizable which shows $M=0$.
We deduce that $\dim W' \leq 1$.
\item Any matrix $M \in W$ is diagonalizable with $M=\begin{bmatrix}
K(M) & C(M) \\
0 & \alpha(M)
\end{bmatrix}$ hence $K(W)$ is a diagonalizable subspace of $\Mat_{n-1}(\K)$, therefore
$\dim K(W) \leq \dbinom{n}{2}$.
\end{itemize}
However $\dbinom{n}{2}+(n-1)+1=\dbinom{n+1}{2}$, and we have assumed that $\dim V=\dbinom{n+1}{2}$.
It follows that:
$$\dim K(W)=\dbinom{n}{2} \quad ; \quad \dim W'=1$$
and $\alpha$ is an isomorphism from $W'$ to $\K$.
Therefore we find a column $C_1 \in \Mat_{n-1,1}(\K)$ such that $W'$ is generated by
$\begin{bmatrix}
0 & C_1 \\
0 & 1
\end{bmatrix}$.
Replacing $V$ with $P^{-1}VP$ for $P=\begin{bmatrix}
I_{n-1} & C_1 \\
0 & 1
\end{bmatrix}$, we are reduced to the situation where $E_{n,n} \in W'$ hence $W'=\Vect(E_{n,n})$.

It follows that $K(W)$ is a maximal diagonalizable subspace of $\Mat_{n-1}(\K)$.
The induction hypothesis thus yields:
\begin{itemize}
\item that $\K$ is formally real and Pythagorean;
\item that there exists a non-singular matrix $Q \in \GL_{n-1}(\K)$ such that $Q\,K(W)\,Q^{-1}=\calS_{n-1}(\K)$.
\end{itemize}
We then find that $P\,V\,P^{-1}$ satisfies all the assumptions from Proposition \ref{finalisation}
hence is conjugate to $\calS_n(\K)$, which shows that $V$ is itself conjugate to $\calS_n(\K)$.
This proves Theorem \ref{structmaxsub} by induction.

\subsection{On non-maximal diagonalizable subspaces}\label{last}

In this short section, we wish to warn the reader that not every diagonalizable subspace of
$\Mat_n(\R)$ may be conjugate to a subspace of $\calS_n(\R)$.
Consider indeed the linear subspace $V$ generated by the matrices
$$A:=\begin{bmatrix}
0 & 0 & 0 \\
0 & 1 & 0 \\
0 & 0 & -1
\end{bmatrix} \quad \text{and} \quad B:=\begin{bmatrix}
0 & 1 & 0 \\
0 & 0 & 1 \\
0 & 1 & 0
\end{bmatrix}.$$
A straightforward computation shows that $a\cdot A+b\cdot B$ has three distinct eigenvalues in $\R$
(namely $0$ and $\pm \sqrt{a^2+b^2}$) for every $(a,b)\in \R^2 \setminus \{(0,0)\}$. Therefore
$V$ is a diagonalizable subspace of $\Mat_3(\R)$. Assume that $V$ is conjugate to a subspace of $\calS_3(\R)$:
then there would be a definite positive symmetric bilinear form $b$ on $\R^3$ such that
$X \mapsto AX$ is self-adjoint. The eigenspaces of $A$ would then be mutually $b$-orthogonal, and the same would hold for $B$.
Denote by $(e_1,e_2,e_3)$ the canonical basis of $\R^3$. Then we would have $\{e_1\}^{\bot_b}=\Vect(e_2,e_3)$.
However $\Vect(e_1)$ is also an eigenspace for $B$ therefore the other two eigenspaces of $B$ should be included in
$\{e_1\}^{\bot_b}=\Vect(e_2,e_3)$, hence $\Vect(e_2,e_3)$ should be their sum. However this fails because $\Vect(e_2,e_3)$
is not stabilized by $B$. This \emph{reductio ad absurdum} shows that $V$ is not conjugate to any subspace of $\calS_3(\R)$.
In particular, there are diagonalizable subspaces of $\Mat_3(\R)$ which are maximal for the inclusion of diagonalizable
subspaces but not in the sense of this paper.

\section{Preserving real diagonalizable matrices}\label{linpres}

Recall that two matrices $A$ and $B$ of $\Mat_n(\K)$ are said to be simultaneously diagonalizable if
there exists a non-singular $P \in \GL_n(\K)$ such that both $PAP^{-1}$ and $PBP^{-1}$ are diagonal.
We will use a recent theorem of Bogdanov and Guterman \cite{Bogdanov}:

\begin{theo}[Bogdanov, Guterman]\label{bogdanov}
Let $f$ be an automorphism of the vector space $\Mat_n(\K)$ and assume that for every
pair $(A,B)$ of simultaneously diagonalizable matrices, $(f(A),f(B))$ is a pair of simultaneously diagonalizable matrices.
Then there exists a non-singular matrix $P \in \GL_n(\K)$, a linear form $\lambda$ on $\Mat_n(\K)$ and a scalar $\mu \in \K^*$
such that $\lambda(I_n) \neq -\mu$, and
$f=\varphi_{\lambda,P,\mu}$ or $f=\psi_{\lambda,P,\mu}$.
\end{theo}

In order to use this, we will see that, under the assumptions of Theorem \ref{linprestheo},
any non-singular linear preserver of diagonalizability must also preserve simultaneous diagonalizability.
This will involve Theorem \ref{structmaxsub} and the following generalization of the
singular value decomposition:

\begin{lemme}[Singular value decomposition theorem (SVD)]\label{SVD}
Let $n \geq 2$. Assume that $\calS_n(\K)$ is a diagonalizable subspace of $\Mat_n(\K)$. \\
Then, for every $P \in \GL_n(\K)$, there exists a triple $(O,U,D) \in O_n(\K)^2 \times D_n(\K)$ such that $P=ODU$.
\end{lemme}

\begin{proof}
The proof is quite similar to the standard one (i.e.\ the case $\K=\R$). \\
Theorem \ref{structmaxsub} shows that $\K$ is formally real and Pythagorean. \\
Let $P \in \GL_n(\K)$. Notice that all the eigenvalues of the non-singular symmetric matrix
$P^TP$ are squares. Let indeed $X$ be an eigenvector of $P^TP$ and $\lambda$ its corresponding eigenvalue.
Then $(PX)^T(PX)=\lambda\,X^TX$. However, both $X^TX$ and $(PX)^T(PX)$ are sums of squares, hence squares,
and $X^TX \neq 0$ since $X \neq 0$ and $\K$ is formally real. It follows that
$\lambda=\frac{(PX)^T(PX)}{X^TX}$ is a square.

Next, the bilinear form
$(X,Y) \mapsto X^T (P^TP)Y=(PX)^T(PY)$ is clearly symmetric hence
$X \mapsto (P^TP)X$ is a self-adjoint operator for the quadratic form
$q : X \mapsto X^TX$. It follows that the eigenspaces of $P^TP$ are mutually $q$-orthogonal. However, since
$\K$ is formally real and Pythagorean, the values of $q$ on $\K^n \setminus \{0\}$ are non-zero squares
hence each eigenspace of $P^TP$ has a $q$-orthonormal basis.

Since the assumptions show that $P^TP$ is diagonalizable,
we deduce that there is an orthogonal matrix $O_1$ and non-zero scalars $\lambda_1,\dots,\lambda_n$ such that $P^TP=O_1D_1O_1^{-1}$ with $D_1:=\Diag(\lambda_1^2,\dots,\lambda_n^2)$. Setting
$S:=O_1DO_1^{-1}$, where $D:=\Diag(\lambda_1,\dots,\lambda_n)$, we find that $S$ is non-singular, symmetric and $S^2=P^TP$.
It is then easily checked that $O_2:=PS^{-1}$ is orthogonal. The matrices $O:=O_2O_1$ and $U:=O_1^{-1}$
are then orthogonal and one has $P=ODU$.
\end{proof}

With the SVD, we may now prove the following result :

\begin{prop}\label{minintersection}
Let $n \geq 2$. Assume that $\calS_n(\K)$ is a diagonalizable subspace of $\Mat_n(\K)$.
Let $V$ and $W$ be two maximal diagonalizable subspaces of $\Mat_n(\K)$.
Then $\dim (V \cap W) \geq n$. If $\dim (V \cap W)=n$, then $V \cap W$ is conjugate to
$D_n(\K)$.
\end{prop}

\begin{proof}
By Theorem \ref{structmaxsub}, we lose no generality assuming that $V=\calS_n(\K)$,
in which case we know that $W=P\,\calS_n(\K) P^{-1}$ for some $P \in \GL_n(\K)$,
and we may then find a triple $(O,U,D) \in O_n(\K)^2 \times D_n(\K)$ such that $P=ODU$.
It follows that $V=\calS_n(\K)$ and $W=OD\calS_n(\K)D^{-1}O^{-1}$.
Conjugating both subspaces by $O^{-1}$, we may assume $V=\calS_n(\K)$ and $W=D\calS_n(\K)D^{-1}$.
In this case, we clearly have $D_n(\K) \subset V \cap W$ and the claimed results follow readily.
\end{proof}

We may now prove Theorem \ref{linprestheo}.

\begin{proof}[Proof of Theorem \ref{linprestheo}]
Let $A$ and $B$ be two simultaneously diagonalizable matrices of $\Mat_n(\K)$.
Replacing $f$ with $M \mapsto f(PMP^{-1})$ for some well-chosen $P \in \GL_n(\K)$,
we lose no generality assuming that $A$ and $B$ are both diagonal.
Consider the diagonal matrix $D:=\Diag(1,2,\dots,n)\in \Mat_n(\K)$ (since $\K$ is formally real, we naturally identify the ring of integers $\Z$
with a subring of $\K$). Then $D$ is non-singular. Set $V:=\calS_n(\K)$ and $W:=D^{-1}\,\calS_n(\K)\,D$.
Then $V \cap W$ consists of all the symmetric matrices $A=(a_{i,j})_{1 \leq i,j \leq n} \in \calS_n(\K)$ such that
$\frac{i}{j}\,\,a_{i,j}=\frac{j}{i}\,a_{j,i}$ for every $(i,j)\in \lcro 1,n\rcro^2$,
i.e.\ such that $(i^2-j^2)\,a_{i,j}=0$ for every $(i,j)\in \lcro 1,n\rcro^2$. Since $\K$
is formally real, one has $i \neq -j$ for every $(i,j) \in \lcro 1,n\rcro^2$, and we deduce that
$V \cap W=D_n(\K)$. Since $A$ and $B$ are both diagonal,
the matrices $f(A)$ and $f(B)$ both belong to $f(V\cap W)$. Since $f$ is one-to-one and linear,
one has $f(V\cap W)=f(V)\cap f(W)$ and $\dim f(V \cap W)=\dim (V \cap W)=n$.
However, the assumptions on $f$ show that $f(V)$ and $f(W)$ are both maximal diagonalizable subspaces of $\Mat_n(\K)$.
It thus follows from Proposition \ref{minintersection} that $f(V) \cap f(W)$ is conjugate to $D_n(\K)$,
hence $f(A)$ and $f(B)$ are simultaneously diagonalizable. Theorem \ref{bogdanov} then yields the claimed results.
\end{proof}

Let us finish with a strong linear preserver theorem:

\begin{theo}\label{stronglinprestheo}
Let $n \geq 2$ be an integer and assume that $\calS_n(\K)$ is diagonalizable.
Let $f$ be an endomorphism of the vector space $\Mat_n(\K)$ such that, for every
$M \in \Mat_n(\K)$, the matrix $f(M)$ is diagonalizable if and only if $M$ is diagonalizable.
Then there exists a non-singular $P \in \GL_n(\K)$, a linear form $\lambda$ on $\Mat_n(\K)$ and a non-zero scalar $\mu$ such that
$f=\varphi_{\lambda,P,\mu}$ or $f=\psi_{\lambda,P,\mu}$.
\end{theo}

Notice conversely that $\varphi_{\lambda,P,\mu}$ and $\psi_{\lambda,P,\mu}$
satisfy the previous assumptions for any $P \in \GL_n(\K)$, any linear form $\lambda$ on $\Mat_n(\K)$ and any non-zero scalar $\mu \in \K^*$.

\begin{proof}
We will reduce the situation to the case $f$ is one-to-one, and the result will then follow directly from Theorem \ref{linprestheo}. \\
We lose no generality assuming $I_n \not\in \Ker f$. If indeed $f(I_n)=0$, then
we choose a linear form $\delta$ on $\Mat_n(\K)$ such that $\delta(I_n) \neq 0$;
the map $g : M \mapsto f(M)+\delta(M)\cdot I_n$ then satisfies the same assumptions as $f$
and the conclusion must hold for $f$ if it holds for $g$.
Assume then that $I_n\not\in \Ker f$, and assume furthermore that $\Ker f$ contains a non-zero matrix
$A$. Pre-composing $f$ with an appropriate conjugation, we then lose no generality assuming that $A=\Diag(\lambda_1,\dots,\lambda_n)$ for some
list $(\lambda_1,\dots,\lambda_n)\in \K^n$ with $\lambda_1 \neq \lambda_2$.
Then, for every diagonalizable matrix $B$, the matrix $f(B)=f(A+B)$ is diagonalizable hence $A+B$ is also diagonalizable.
Set $B:=\begin{bmatrix}
M & 0 \\
0 & N
\end{bmatrix}$ with $M:=\begin{bmatrix}
-\lambda_1 & 1 \\
0 & -\lambda_2
\end{bmatrix}$ and $N:=\Diag(-\lambda_3,\dots,-\lambda_n)$, and note that $B$ is diagonalizable whereas $A+B$ is
not since it is nilpotent and non-zero. This contradiction shows that $\Ker f=\{0\}$, and Theorem \ref{linprestheo} then yields the
desired conclusion.
\end{proof}

Still open is the question of the determination of the endomorphisms
$f$ of $\Mat_n(\K)$ such that $f$ stabilizes $\DG_n(\K)$:
clearly, the $\varphi_{\lambda,P,\mu}$ and $\psi_{\lambda,P,\mu}$'s always qualify; also,
picking an arbitrary diagonalizable subspace $V$ of $\Mat_n(\K)$, any linear map $f : \Mat_n(\K) \rightarrow V$ also qualifies.
We do not know whether all the solutions have one of the above forms.

\begin{Rem}[The case $n=2$]
Let $\K$ be an arbitrary field and assume that not every matrix of $\calS_2(\K)$ is diagonalizable.
Hence Theorem \ref{structmaxsub} shows that every diagonalizable subspace of $\Mat_2(\K)$
has dimension at most $2$. Let $V$ be such a subspace. If $\dim V=2$, then
$V+\Vect(I_2)$ is still diagonalizable hence $I_2 \in V$. Completing $I_2$ into a basis $(I_2,A)$ of $V$,
we easily see, using the fact that $A$ is diagonalizable that $V$ is conjugate to $D_2(\K)$.
It easily follows that any linear automorphism $f$ of $\Mat_2(\K)$ which preserves diagonalizability must also
preserve simultaneous diagonalizability.
We deduce that Theorem \ref{linprestheo} actually holds for an arbitrary field if $n=2$.
\end{Rem}

\end{document}